\documentclass[reqno,10pt]{amsart}
\usepackage{graphicx}

\newcommand \RR         {\mathbb{R}}
\newcommand \del        \partial
\newcommand \eps        \epsilon
\newcommand \lam        \lambda
\newcommand \be         {\begin{equation}}
\newcommand \ee         {\end{equation}}

\newcommand \CC      {\mathcal C}

\newcommand\GG   { G}
\newcommand\ld    {\lambda}

\newcommand\ZZ    {\mathbf Z}
\newcommand\FF   {{\mathcal F}}
\newcommand\UU    {{\mathcal U}}

\newcommand\pt   {\partial}

\newcommand{\rf}[1]{(\ref{#1})}

\numberwithin{equation}{section}
\newtheorem{theorem}{Theorem}[section]
\newtheorem{lemma}[theorem]{Lemma}
\newtheorem{proposition}[theorem]{Proposition}

\newtheorem{remark}[theorem]{Remark}

\begin{document}

\title[The minimum entropy principle for compressible fluid flows]
{\small The minimum entropy principle for compressible fluid flows in a nozzle
with discontinuous cross-section\footnote{\tiny{To appear in: Mathematical Modelling and Numerical Analysis (M2AN).}}}
\author[D. K{\tiny r\"oner}, P.G. L{\tiny e}F{\tiny loch,}  and M.D. T{\tiny hanh}]
{\bf
D{\tiny ietmar} K{\tiny r\"oner,}  P{\tiny hilippe} G. L{\tiny e}F{\tiny loch,} {\tiny and}
  M{\tiny ai}-D{\tiny uc} T{\tiny hanh}}

\address{Dietmar Kr\"oner,
Institute of Applied Mathematics, University of Freiburg,
        Hermann-Herder Str.~10, 79104 Freiburg, Germany.}
        \email{Dietmar.Kroener@mathematik.uni-freiburg.de.}

\address{ Philippe G. L{\small e}Floch,
Laboratoire Jacques-Louis Lions \& Centre National de la Recherche
Scientifique, Universit\'e de Paris VI, 4 Place Jussieu, 75252 Paris,
France.} \email{LeFloch@ann.jussieu.fr}

\address{Mai Duc Thanh,
Department of Mathematics, International University, Quarter 6,
Linh Trung Ward, Thu Duc District, Ho Chi Minh City, Vietnam}
\email{MDThanh@hcmiu.edu.vn}


\begin{abstract} We consider the Euler equations for compressible fluids
in a nozzle whose cross-section is variable and may contain discontinuities.
We view these equations as a hyperbolic system in nonconservative form
and investigate weak solutions in the sense of Dal Maso, LeFloch, and Murat.
Observing that the entropy equality has a fully conservative form,
we derive a minimum entropy principle satisfied by entropy solutions.
We then establish the stability of a class of numerical approximations for this system.
\end{abstract}
\maketitle

\section{Introduction}
\label{IN}

Compressible flows in a nozzle with variable cross-section $a=a(x)$ are described by
the Euler equations, which under symmetry assumptions take the following form \cite{CourantFriedrichsbook}
\be
\aligned &\del_t(a\rho) + \del_x(a\rho u) = 0,
\\
&\del_t(a\rho u) + \del_x(a \rho u^2) + a \, \del_x p = 0,
\\
&\del_t (a\rho e) + \del_x(a u (\rho e + p))  =  0,
\quad x\in\RR, \, t>0,
\endaligned
\label{1.1}
\ee
where the main unknowns are the fluid velocity $u$ and the thermodynamic variables $\eps$ and
$\rho$. Here,  $e = \eps + u^2/2$ is the total energy,
and we also define the specific volume $v = 1/\rho$ and the pressure $p$,
the temperature $T$, and the specific entropy $S$. These variables
are related via the equation of state of the fluid under consideration, for instance in the form
$$
p=p(\rho,S).
$$
Provided $p$ is a monotone increasing function of $\rho$ (for fixed $S$), the system above
is a hyperbolic system of balance law with variable coefficients.

To deal with discontinuous cross section, following
\cite{LeFloch89} we supplement \eqref{1.1} with the ``trivial''
equation \be \del_t a = 0, \label{1.2} \ee so that the whole set
of equations can be written as a hyperbolic system in
nonconservative form \be \del_t U + A(U) \, \del_x U = 0.
\label{1.3} \ee In consequence, at least within the regime where
the system is strictly hyperbolic, the theory of such systems
developed by Dal Maso, LeFloch, and Murat
\cite{DalMasoLeFlochMurat} (and also \cite{LeFlochLiu,
LeFlochbook, LeFloch04}) applies, and provide the existence of
entropy solutions to the Riemann problem (a single discontinuity
separating two constant states as an initial data), as well as to
the Cauchy problem (for solution with sufficiently small total
variation). More recently, LeFloch and Thanh \cite{LeFlochThanh,LeFlochThanh2}
solved the Riemann problem for arbitrary data, including the
regime where the system fails to be globally strict hyperbolicity (i.e., the
resonant case).  The Riemann problem was also solved by a different approach by Andrianov and Warnecke
\cite{AndrianovWarnecke}. For earlier work on resonant systems,
see also \cite{MarchesinPaes-Leme,IsaacsonTemple95,
IsaacsonTemple92,GoatinLeFloch}.

In the present paper, we pursue the analysis of the Euler equations in a general nozzle,
and establish several properties of solutions.
We derive the entropy inequality that must be satisfied by weak solutions.
The entropy inequality is found to have a fully conservative form, so that the
notion in nonconservative product is not needed to state the entropy inequality.
In turn, we can obtain a generalization of the so-called minimum
entropy principle, originally established by Tadmor \cite{Tadmor86} (see also \cite{Tadmor84})
for plane symmetric fluid flows.  In particular, it follows that the specific entropy is non-increasing in time.

These properties are important as far as the stability of numerical schemes is concerned, and
one of our main results is a proof that a scheme proposed by Kr\"oner and Thanh \cite{KroenerThanh2}
satisfies a variant of the minimum entropy principle.
Recall that the discretization of systems of balance laws is particularly delicate, and was addressed
by many authors \cite{GreenbergLeroux,
BotchorishviliPerthameVasseur, BotchorishviliPironneau, Gosse00, Bouchut,
AudusseBouchutBristeauKleinPerthame}.
We show that the scheme under consideration
not only preserve equilibrium states, but also preserves the positivity of the density and
the minimum entropy principle.

%------------------------------------------------------------------------------------------------------------------

\section{Entropy inequality for nozzle flows}
\label{section2}

Consider the system \rf{1.1} supplemented with a given equation of state
$p=p(\rho, \eps)$,
and set $U = (\rho, \rho u, \rho e)$. The flux function and the right-hand side in \rf{1.1}
can be expressed as functions of $U$:
$$
\aligned
(\rho u, \rho u^2 + p(\rho, \eps), u (\rho e + p(\rho, \eps))) &:= f(U),\\
 (0,p(\rho, \eps),0) &:= g(U).\\
 \endaligned
$$
The system under consideration can therefore be written as
$$
\del_t (aU) +  \del_x (af(U))  = g(U)\frac{da}{dx}.
$$

Recall that weak solutions to this nonconservative systems are defined in the sense of
Dal~Maso, LeFloch, and Murat \cite{DalMasoLeFlochMurat}. On the other hand,
as we will see later on, the (mathematical) entropy inequality associated with this system has
a conservative form and does make sense in the framework of distributions.

Consider first the one-dimensional gas dynamics equations corresponding to a constant
function $a$ ($x\in\RR$, $t>0$)
\be
\aligned &\del_t\rho + \del_x(\rho u) \, = \, 0,
\\
&\del_t(\rho u) + \del_x(\rho u^2 + p) \, = \, 0,
\\
&\del_t (\rho e) + \del_x(u (\rho e + p)) \, = \, 0.
\endaligned
\label{2.1}
\ee
On one hand, the notion of entropy is motivated from physics, and the physical entropy is
$\UU = \rho S$, where $S$ is the specific entropy.
On the other hand, as was shown by Harten et al \cite{HartenLaxLevermoreMorokoff98},
necessary and sufficient conditions for a
twice differentiable function $\UU_c$ of  the form
\be
\UU_c = \rho g(S),
\label{2.2}
\ee
to be an entropy of the usual gas dynamics equations is that $g(S)$ satisfies the following properties:
\begin{itemize}
\item[(i)]
$g(S)$ is strictly decreasing as  function of $S$;
\item[(ii)]
$g(S)$ is strictly convex as function of $(1/\rho, \eps)$.
\end{itemize}
Moreover, the system under consideration is strictly hyperbolic if and only if
it admits an entropy of the form \rf{2.2}.

Consider the Navier-Stokes equations describing a viscous fluid flow
in a nozzle with smooth area function $a_\nu=a_\nu(x)$
\be
\aligned
&\del_t(a_\nu\rho_\nu) + \del_x(a_\nu\rho_\nu u_\nu) \, = \, 0,
\\
&\del_t(a_\nu\rho_\nu u_\nu) + \del_x(a_\nu(\rho_\nu u_\nu^2 + p_\nu)) \, = \, p_\nu\del_x a_\nu
+ \nu \,
\del_x (b_\nu  \del_x u_\nu),
\\
&\del_t (a_\nu\rho_\nu e_\nu) + \del_x(a_\nu u_\nu (\rho_\nu e_\nu + p_\nu)) \,
= \nu \, \del_x(b_\nu u_\nu\del_x u_\nu),
\endaligned
\label{2.3}
\ee
where $b_\nu =b_\nu(x)\ge 0$ is given and
$\nu$ denotes the viscosity coefficient.

We consider the limit as $\nu$ tends to zero.
For simplicity, we drop the subscript $\nu$ and derive the equation for the specific entropy.
To this end, we assume that the internal energy is given by an equation of state $\eps = \eps(\rho,S)$.
On one hand, thanks to the equation of conservation of mass of \rf{2.3},
the equation of momentum in \rf{2.3} can be written as
\be
a\rho(u_t  +  uu_x) + p_xa - \nu \, (b\,  u_x)_x = 0.
\label{2.4}
\ee
On the other hand, the equation of energy in \rf{2.3} can be written as
\be
a\rho e_t + e\big( (a\rho)_t  + (a\rho u)_x\big) + a\rho ue_x + (aup)_x \, = \nu \, (b \, u u_x))_x.
\label{2.5}
\ee
The second term on the left-hand side of \rf{2.5} is equal to zero due to the conservation of mass.
Using the thermodynamical identity
$d\eps = TdS - pdv$ and $v = 1/\rho$,
we can re-write the equation \rf{2.5} as
$$
a\rho T(S_t + uS_x) + \frac{ap}{\rho}(\rho_t + u\rho_x) + (aup)_x + a\rho u(u_t + uu_x)
= \nu \, (b_\nu \, u u_x ))_x.
$$
Or, after arranging terms, we obtain the equation of energy
\be
\aligned
& a\rho T(S_t + uS_x) + \frac{p}{\rho}\big( (a\rho)_t + (au\rho)_x\big)
+ u\big(a\rho (u_t + uu_x)
+ ap_x - \nu \, (b\, u_x)_x\big)  \\
& = \nu \, b\, u_x^2.
\endaligned
\label{2.7}
\ee
The second and the third term on the left-hand side of \rf{2.7} are equal to zero
by the conservation of mass of \rf{2.3} and the momentum equation \rf{2.4}.
Thus, we deduce from \rf{2.7} that
the specific entropy $S$ should satisfy
$$
\del_t S + u \del_x S = \frac{b\nu }{a\rho T} u_x^2.
$$
Let $g(S)$ be any smooth function of $S$ satisfying $g'(S) \le 0$.
Multiplying the above equation by $a\rho g'(S)$, we obtain
$$
a\rho\del_t g(S) + a\rho u \del_x g(S) =
\nu \,
a\rho g'(S)\Big(\frac{b}{a \rho T} u_x^2 \Big).
$$
Multiplying the conservation of mass of \rf{2.3}  by $g(S)$ and then summing up
with the above equation, we find
$$
\del_t (a\rho g(S)) + \del_x (a\rho u g(S))
=
\nu \,
\rho g'(S)\Big(\frac{b}{\rho T} u_x^2 \Big).
$$

It is clear that
if the system  \rf{2.3} admits a sequence of smooth solution $U_\nu$ uniformly bounded in amplitude
and converging
almost everywhere to a limit $U=(\rho, \rho u, \rho e)$ when $\nu$ tends to zero,
then, the function $U$ satisfies the {\sl entropy inequality}
\be
\del_t (a\rho g(S)) + \del_x (a\rho u g(S)) \le 0,
\label{2.8}
\ee
where $g$ is any function satisfying the above items (i)-(ii).
Recall that a weak solution of \rf{1.1} satisfying the entropy inequality \rf{2.8} in the distributional
sense is called an entropy solution.

%--------------------------------------------------------------------------------------------------------------

\section{Minimum entropy principle for nozzle flows}

We always assume that the fluid is in local thermodynamic equilibrium, so that:
$$
\text{The function $(v,\eps)\mapsto \eps(v, S)$, is strictly convex.}
$$
This assumption is equivalent to the requirement that the function $(v,\eps)\mapsto S(v,\eps)$ is strictly concave.

Thanks to the divergence form of the entropy inequality \rf{2.8}, we can establish a minimum entropy principle.
We begin with the entropy inequality in a generalized form and check
that the entropy inequality for \rf{1.1}-\rf{1.2} coincides with \rf{2.8} for the entropy pair
\be
(\UU,\FF)  = (a\rho g(S),a\rho ug(S)),
\label{3.1}
\ee
where the functions $g$ satisfy the assumptions (i)-(ii) in  Section~\ref{section2}.

Consider the hyperbolic system in nonconservative form
\be
\del_t U + A(U)\,  \del_x U = 0.
\label{3.2}
\ee
The entropy inequality for \rf{3.2} has the form
$$
\del_t \UU (U) + \Big[D_U\UU (U) \,  A(U(.,t)) \del_x
U(.,t)\Big]_\phi  \le 0,
$$
where $\phi$ is a given Lipschitz family of paths, and $\UU$ is a convex function satisfying
$$
D_U^2\UU (U) \,  A(U) = A(U)^T\,  D_U^2\UU (U).
$$
Basic properties of the nonconservative product imply that if there exists a function $\FF$ such that
\be
D_U\UU\,  A(U) = D_U\FF(U),
\label{3.3}
\ee
then the nonconservative product
$\Big[D_U\UU (U) \,  A(U(.,t)) \del_x U(.,t)\Big]_\phi$
reduces to the usual one in divergence form, and is independent of the path $\phi$.
Consequently, the entropy inequality takes the divergence form
\be
\del_t \UU (U) + \del_x \FF (U)  \le 0
\label{3.4}
\ee
in the sense of distributions.

We now check that the entropy inequality for
\rf{1.1} can be reduced to the divergence form \rf{3.4}
for all entropy pairs of the form \rf{3.1}. This will establish \rf{3.3}.

The system \rf{1.1}-\rf{1.2} can be written in the nonconservative
form \rf{3.2}, with $U = (a\rho, a\rho u, a\rho e, a) := (w_1,w_2,w_3,w_4)$
Replacing the expression of $U$ in the system \rf{1.1}-\rf{1.2}, we obtain
$$
\aligned &\del_t w_1 + \del_x w_2 \, = \, 0,
\\
&\del_t w_2 + \dfrac{w_2}{w_1}\del_x w_2 + w_2 \del_x\big(\dfrac{w_2}{w_1}\big) + w_4\del_xp \, = \, 0,
\\
&\del_t w_3 + w_2 \del_x\big(\dfrac{w_3}{w_1}\big) +\dfrac{w_3}{w_1}\del_x w_2 + \del_x(aup) \, = \, 0,\\
& \del_t w_4  = 0.
\endaligned
$$
After a tedious but straightforward calculation, we arrive at the following system
$$
\aligned
&\del_t w_1 + \del_x w_2 \, = \, 0,
\\
&\del_t w_2 +
\Big(p_\rho - \frac{p_\eps w_3w_4}{w_1^2} + \frac{p_\eps w_2^2w_4}{w_1^3} - \frac{w_2^2}{w_1^2}\Big)\,  \del_xw_1
+\Big(\frac{2w_2}{w_1} - \frac{p_\eps w_2w_4}{w_1^2}\Big)\del_x w_2 \\
&\qquad\qquad + \frac{p_\eps w_4}{w_1}\del_x w_3
 - \frac{p_\rho w_1}{w_4}\del_x w_4  \, = \, 0,\\
&\del_t w_3 +
    \Big(     \frac{w_2w_4}{w_1}\Big(\frac{p_\rho}{w_4} - \frac{p_\eps w_3}{w_1^2}
     + \frac{p_\eps w_2^2}{w_1^3} - \frac{p}{w_1}\Big) - \frac{w_2w_3}{w_1^2}\Big)
     \,  \del_xw_1\\
&\qquad +\Big(\frac{w_3 + pw_4}{w_1} - \frac{p_\eps w_2^2w_4}{w_1^3}\Big)\del_x w_2
+ \Big(\frac{p_\eps w_2w_4}{w_1^2} + \frac{p+w_2}{w_1}\Big)\del_x w_3\\
&\qquad\qquad +\frac{w_2}{w_1}\Big(p -  \frac{p_\rho w_1}{w_4}\Big)\del_x w_4  \, = \, 0,\\
& \del_t w_4  = 0,\quad x\in\RR, \, t>0.\\
\endaligned
$$
This system has the canonical form \rf{3.2}, where
the matrix $A(U) = \big(a_{ij}(U)\big)$
is given by:
\be
\aligned
&a_{11} = 0,\quad a_{12} = 1,\quad a_{13} = 0,\quad a_{14} = 0,\\
&a_{21} = p_\rho - \frac{p_\eps w_3w_4}{w_1^2} + \frac{p_\eps w_2^2w_4}{w_1^3} - \frac{w_2^2}{w_1^2}
= p_\rho - u^2 + \frac{p_\eps}{\rho}(u^2-e),\\
&a_{22} = \frac{2w_2}{w_1} - \frac{p_\eps w_2w_4}{w_1^2} = 2u - \frac{p_\eps u}{\rho},\\
&a_{23} = \frac{p_\eps w_4}{w_1} = \frac{p_\eps}{\rho},\\
&a_{24} = - \frac{p_\rho w_1}{w_4} = -p_\rho\rho,\\
&a_{31} = \frac{w_2w_4}{w_1}\Big(\frac{p_\rho}{w_4} - \frac{p_\eps w_3}{w_1^2}
     + \frac{p_\eps w_2^2}{w_1^3} - \frac{p}{w_1}\Big) - \frac{w_2w_3}{w_1^2}
     = u(p_\rho - e + \frac{p_\eps (u^2 - e) - p}{\rho}),\\
&a_{32} = \frac{w_3 + pw_4}{w_1} - \frac{p_\eps w_2^2w_4}{w_1^3} = e + \frac{p}{\rho} - \frac{p_\eps u^2}{\rho},\\
&a_{33} = \frac{p_\eps w_2w_4}{w_1^2} + \frac{p+w_2}{w_1} = \frac{p_\eps u}{\rho} + u,\\
&a_{34} = \frac{w_2}{w_1}\Big(p -  \frac{p_\rho w_1}{w_4}\Big) = u(p-p_\rho\rho ),\\
&a_{41} = a_{42} = a_{43} =a_{44} = 0.\\
\endaligned
\label{3.6}
\ee
These coefficients will be needed in the proof of the forthcoming theorem.

\begin{proposition} \label{theo41}
Consider the system \rf{1.1}-\rf{1.2} in the form  \rf{3.2}, \rf{3.6}.
Let $g$ be any function satisfying the hypotheses i) and ii) in Section~3.
Then the function $\UU  = a\rho g(S)$ of the conservative variables $(a\rho, a\rho u, a\rho e, a)$
is convex. Moreover, it satisfies
$$
D_U\UU\,  A(U) = D_U\FF(U),\quad \FF(U) = a\rho ug(S),
$$
which implies that $(\UU, \FF)$ is an entropy-pair of the Euler system. Consequently, the
entropy inequality in the sense of nonconservative products can be written in the
 divergence form
\be
(a\rho g(S))_t  + (a\rho ug(S))_x \le 0.
\label{3.8}
\ee
\end{proposition}

\begin{proof}
First, as is shown in \cite{HartenLaxLevermoreMorokoff98}, the function $\rho g(S)$ is convex in the variable
\newline
$(\rho, \rho u, \rho e)$. Therefore, the function $a\rho
g(S)=\bar\rho g(S), \ \bar\rho :=a\rho$, is convex in the variable
$(\bar\rho,\bar\rho u,\bar\rho e)$. Since this function $\bar\rho
g(S)$ can be seen as dependent only on the first three variables
$(\bar\rho,\bar\rho u,\bar\rho e)$, it can thus be seen as
independent of $a$. Therefore, it is also convex in the variable
$(\bar\rho,\bar\rho u,\bar\rho e, a)= (a\rho, a\rho u, a\rho e,
a)$.

Second, the equation of state for the specific entropy being written as
$S = S(\rho, \eps)$,
 a straightforward calculation shows that
\be
\UU(U) = a\rho g(S) = w_1g\Big(S\Big(\frac{w_1}{w_4}, \frac{w_3}{w_1}
 - \frac{1}{2}\big(\frac{w_2}{w_1}\big)^2\Big)\Big),\quad \FF(U) = u\UU(U).
\label{3.10}
\ee
Using the thermodynamic identity
$d\eps = TdS - pdv = TdS + \frac{p}{\rho^2} d\rho$,
we have
$$
S_\eps = \frac{1}{T},\quad S_\rho = \frac{-p}{T\rho^2}.
$$
Note also that
$$
{\mathcal U}_{w_3} = {g'(S) \over T}, \qquad
{\mathcal U}_{w_4} = p \, {g'(S) \over T}.
$$
Therefore, it follows that
\be
\aligned
D_U\UU(U) &= \Big(g(S) + w_1g'(S)\big(\frac{S_\rho}{w_4} + \frac{S_\eps}{w_1^2}(\frac{w_2}{2} - w_3)\big),\\
&\qquad\qquad -g'(S)S_\eps\frac{1}{2}, g'(S)S_\eps, -g'(S)S_\rho\frac{w_1^2}{w_4^2}\Big)\\
& = \big(g(S) + \frac{g'(S)}{T}(-\frac{p}{\rho} - e + u^2), -\frac{g'(S)}{2T}, \frac{g'(S)}{T}, \frac{g'(S)p}{T}\big).\\
\endaligned
\label{3.11}
\ee
and
\be
\aligned
D_U\FF(U) &= uD_U\UU(U) + \UU(U) D_U\big(\frac{w_2}{w_1}\big)\\
& = (u\UU_{w_1} - \frac{w_2}{w_1^2}\UU(U), u\UU_{w_2} + \frac{\UU(U)}{w_1}, u\UU_{w_3}, u\UU_{w_4}).
\endaligned
\label{3.12}
\ee
From \rf{3.6},  \rf{3.11}, and \rf{3.12}, we
claim that \be B := D_U\UU\,  A(U) - D_U\FF(U) = 0. \label{3.13}
\ee Actually,  setting $B = (b_1, b_2, b_3,b_4)$, we have
$$
\aligned
b_1 &= -w_1g'(S)S_\eps\frac{w_2}{w_1} a_{21} + w_1g'(S)\frac{S_\eps}{w_1} a_{31} - u \UU_{w_1} + \frac{w_2}{w_1^2}\UU\\
&= g'(S)S_\eps (a_{31} - \frac{w_2}{w_1}a_{21}) - u \UU_{w_1} + \frac{w_2}{w_1^2}\UU\\
&=g'(S)S_\eps (\frac{-pw_2w_4}{w_1^2} + \frac{w_2^3}{w_1^3} - \frac{w_2w_3}{w_1^2}) - u \UU_{w_1} + \frac{w_2}{w_1^2}\UU,
\endaligned
$$
thus
$$
\aligned
b_1
&=g'(S)S_\eps (u^3 - \frac{pu}{\rho} -ue) - u \UU_{w_1} + \frac{w_2}{w_1^2}\UU\\
&=\frac{g'(S)u}{T} (u^2 - \frac{p}{\rho} -e) - u (g(S) + \frac{g'(S)}{T}(\frac{-p}{\rho} + u^2 - e)) + ug(S)\\
&=0,\\
\endaligned
$$
and
$$
\aligned
b_2 & = g(S) + \frac{g'(S)}{T}(\frac{-p}{\rho} + u^2 - e) -\frac{g'(S)u}{T}(2u - p_\eps\frac{u}{\rho}) +
\frac{g'(S)}{T}(e + \frac{p}{\rho} - p_\eps\frac{u^2}{\rho})\\
&\qquad\qquad - (u\UU_{w_2} + \frac{\UU}{w_1})\\
&= g(S) - \frac{g'(S)u^2}{T} - u\UU_{w_2} - \frac{\UU}{w_1})
= 0,\\
\endaligned
$$
and finally
$$
\aligned
b_3 &=-\frac{g'(S)u}{T}\frac{p_\eps}{\rho} + \frac{g'(S)}{T}(\frac{p_\eps u}{\rho} + u) - u\UU_{w_3} \\
b_4 & = p_\rho\rho{g'(S)u\over T} + u(p-p_\rho\rho) {g'(S)\over T} - u\UU_{w_4} =0.\\
\endaligned
$$
From the above relations we easily check \rf{3.13}, which completes the proof of the proposition.
\end{proof}

%---------------------------------------------------------------------------------------------------------

We are now in a position to establish the minimum entropy principle for gas flows in a nozzle.

\begin{theorem} \label{theo42}
If $U$ is a bounded entropy solution to the system \rf{1.3}, then
it satisfies the {\rm minimum entropy principle}:
$$
\underset{|x|\le R}{\text{\rm inf }} S(x,t)\ge
\underset{|x|\le R+t||u||_{L^\infty}}{\text{\rm inf }} S(x,0).
$$
\end{theorem}

We will need:

\begin{lemma}\label{lem34}
Given a real $p> 1$,
consider the function $g(S):=(S_0-S)^p$,
where $S_0$ is a constant such that $S_0-S>0$ for all $S$
in the domain under consideration. Then:
\begin{itemize}
\item[(i)] $g(S)$  is strictly decreasing and strictly convex as a function of $S$,
\item[(ii)] $g(S)$ is strictly convex as a function of $(v,\eps)$.
\end{itemize}
\end{lemma}

\begin{proof} We have $g'(S)=-p(S_0-S)^{p-1}<0$
and $g''(S)=p(p-1)(S_0-S)^{p-2}>0$, so (i) follow imediately.

Next, since the function $S(v,\eps)$ is strictly concave as a
function of $(v,\eps)$ for $0<\ld<1$, for
$(v_1,\eps_1)\ne (v_2,\eps_2)$ we have
$$
S(\ld(v_1,\eps_1)+(1-\ld)(v_2,\eps_2)) > \ld
S(v_1,\eps_1)+(1-\ld)S(v_2,\eps_2).
$$
Thus, by (i) it follows
$$
\aligned
g\big(S(\ld(v_1,\eps_1)+(1-\ld)(v_2,\eps_2))\big)
& < g\big(\ld
S(v_1,\eps_1)+(1-\ld)S(v_2,\eps_2)\big)
\\
& < \ld g\big(S(v_1,\eps_1)\big)+(1-\ld)g\big(S(v_2,\eps_2)\big),
\endaligned
$$
which establishes (ii).
\end{proof}

\begin{proof}[Proof of Theorem~\ref{theo42}]
Let $g=g(S)$, where $S$ is the specific entropy,
 be any function satisfying the conditions (i) and (ii) stated earlier.
We claim that any bounded entropy solution of the system
\rf{1.1} satisfies
\be
\int_{|x|\le R} \rho (x,t) g(S(x,t)) dx \,
\le \, \int_{|x|\le R + t||u||_{L^\infty}} \rho (x,0) g(S(x,0))
dx, \label{3.18}
\ee
for any nonnegative decreasing function
$g(S)$.

We follow the arguments in the proof of Lemma 3.1 of \cite{Tadmor86}.
We integrate the entropy inequality in the divergence form \rf{3.8} over the truncated cone $\CC = \{ (x,t) |
|x|\le R+(t-\tau)||u||_{L^\infty}, \quad 0\le \tau\le t\}$. Denoting by $(n_x,n_t)$
the unit outer normal of $\CC$, Green's formula yields
$$
\int_{\pt\CC} a\rho g(S) (n_t+un_x)d s \le 0.
$$
The integrals over the top and bottom lines of $\pt \CC$ give the
difference between the left- and the right-hand sides of \rf{3.18}. It follows from the last inequality that
this term is bounded from above by
$$
-\int_{\text{mantle of } \CC} a \, \rho g(S) \, (n_t+un_x) \, d s.
$$
We will show that the last quantity is non-positive. Indeed, on the mantle we have
$$
(n_x,n_t)=(1+||u||_{L^\infty}^2)^{-1/2}(x/|x|, ||u||_{L^\infty}),
$$
thus
$$
n_t+un_x=(1+||u||_{L^\infty}^2)^{-1/2}\Big(||u||_{L^\infty}+
ux/|x|\Big) \ge 0.
$$
By the condition $g(S)>0$, we obtain the desired conclusion.

Now, consider the family of function $g(S) = (S_0-S)^p$, where
$p > 1$ and $S_0$ is a constant satisfying $S+S_0>0$; for instance, $S_0=||S||_{L^\infty}+1$.
To proceed, we use Lemma~\ref{lem34}. The
inequality \rf{3.18} yields
$$
\int_{|x|\le R} \rho (x,t) (S_0-S(x,t))^p dx  \le \int_{|x|\le R +
t||u||_{L^\infty}} \rho (x,0) (S_0-S(x,0))^p dx,
$$
or
$$
\Big(\int_{|x|\le R} \rho (x,t) (S_0-S(x,t))^p dx\Big)^{1/p}  \le
\Big(\int_{|x|\le R + t||u||_{L^\infty}} \rho (x,0) (S_0-S(x,0))^p
dx\Big)^{1/p}.
$$
This means that
\be
\aligned
&||\rho^{1/p} (.,t) (S_0-S(.,t))||_{L^p[-R,R]} \\
&\qquad\qquad \le \,
 ||\rho^{1/p} (.,0) (S_0-S(.,0))||_{L^p[-R - t||u||_{L^\infty}, R + t ||u||_{L^\infty}]}.
 \endaligned
\label{3.19}
\ee

Letting $p\to +\infty$ in the above inequality we obtain
$$
||S_0-S||_{L^\infty[-R,R]} \le ||S_0-S||_{L^\infty[-R -
t||u||_{L^\infty},R + t||u||_{L^\infty}]}.
$$
By definition, we can write this in the form
$$
\underset{|x|\le R}{\text{\rm sup }} (S_0- S(x,t))\le
\underset{|x|\le R+t||u||_{L^\infty}}{\text{\rm sup }}
(S_0-S(x,0)),
$$
or, suppressing the large constant $S_0$,
$$
S_0+\underset{|x|\le R}{\text{\rm sup }} (- S(x,t))\le
S_0+\underset{|x|\le R+t||u||_{L^\infty}}{\text{\rm sup }}
(-S(x,0)).
$$
Eliminating $S_0$ and using (for instance)
$\text{\rm sup } (- S(x,t))=-\text{\rm inf } S(x,t)$, we arrive at the desired result
and the proof of Theorem \ref{theo42} is completed.
\end{proof}

%==================================================================================================================

\section{An entropy stable and well-balanced scheme for fluid flows in a nozzle}

\subsection{Equilibrium states and admissibility criterion}

In this section we investigate various properties of approximate solutions
generated by a finite difference scheme for the Euler equations in a nozzle; this scheme was first proposed
in \cite{KroenerThanh1,KroenerThanh2}. For definiteness and clarity in the presentation,
we consider stiffened gases described by
\be
p = (\gamma - 1) \rho (\eps - \eps_\infty) - \gamma p_\infty, \quad 1 <  \gamma < 5/3,
\label{4.1}
\ee
where $\eps_\infty, p_\infty$ are constants depending on the material under consideration 
with ${p_\infty \ge 0}$.

One key property of the well-balanced scheme under consideration is that it preserves equilibrium states. In this subsection, 
we will recall some basic facts and explain our selection criterion in the construction of
 the right-hand state that can be connected to a given left-hand state by a stationary wave.

Let us observe first that the system \rf{1.1}-\rf{1.2} is non-strictly hyperbolic. More precisely, the phase domain is divided into three sub-domains so that in each of these domains the system is strictly hyperbolic, and along the phase boundary the characteristic fields coincide. We express here all other thermodynamics variables in term of $(\rho, S)$. 
For smooth solutions, the system \rf{1.1}-\rf{1.2} is equivalent to
\be
\aligned
& \rho_t + u\rho_x + \rho u_x + {\rho u\over a} a_x =0,
\\
& u_t+ {p_\rho\over\rho}\rho_x +  uu_x + {p_S\over \rho}S_x =0,
\\
& S_t + uS_x =0,
\\
& a_t = 0.
\endaligned
\label{4.2}
\ee
Thus, in the variable $U = (\rho, u, S,a)$ 
the system \rf{1.1}-\rf{1.2} for smooth flows can be written in the nonconservative form
$U_t + A(U)U_x=0$,
 where
$$
A(U)  =  \left(\begin{matrix} u & \rho & 0 &\dfrac{u \rho}{a}
\\
\dfrac{p_\rho}{\rho} & u  & \dfrac{p_S}{\rho} &0\\
0 & 0  & u &0\\
0&0&0 & 0\\
\end{matrix}\right).
$$
The matrix $A(U)$ admits four real eigenvalues, provided $p_\rho(\rho,S)>0$. Therefore, for $p_\rho(\rho,S)>0$
 the system
\rf{1.1}-\rf{1.2} has four characteristic fields associated with 
the eigenvalues
$$
 \lambda_0 = 0,\quad \lambda_1 = u- \sqrt{p_\rho(\rho,S)}, \quad \lambda_2 = u,\quad \lambda_3 = u+\sqrt{p_\rho(\rho,S)}.
$$
The phase space is decomposed in several regions:
\be
\aligned
&\GG_1 = \{U : \ld_1(U) < \ld_2(U) < \ld_3(U) < \ld_0(U)\},\\
&\GG_2 = \{U : \ld_1(U) < \ld_2(U) < \ld_0(U) < \ld_3(U)\},\\
&\GG_3 = \{U : \ld_1(U) < \ld_0(U) < \ld_2(U) < \ld_3(U)\},\\
&\GG_4 = \{U : \ld_0(U) < \ld_1(U) < \ld_2(U) < \ld_3(U)\},\\
\endaligned
\label{4.3}
\ee
together with isolated surfaces along which the system fails to be strictly hyperbolic:
\be
\aligned
&\Sigma_+ = \{U :  \ld_1(U) = \ld_0(U)\},\\
&\Sigma_0 = \{U :  \ld_2(U) = \ld_0(U)\},\\
&\Sigma_- = \{U :  \ld_3(U) = \ld_0(U)\}.\\
\endaligned
\label{4.4}
\ee
Note that the matrix of the hyperbolic system is diagonalizable on $\Sigma_O$ too.
(This property is relevant for a flow at rest for which $u=0$ even on a discontinuity of $a$.)

Next, we consider some properties of equilibrium states. The entropy is constant across any stationary wave. So, 
we may talk about states by ignoring the component $S$. Suppose that a left-hand state $U_-=(\rho_-,u_-,a_-)$ is given, 
where $a=a_-$ is the value of the cross-section. 
A state $U_+=(\rho_+,u_+,a_+)$ (with the cross-section $a_+$) 
which can be connected with $U_-$ via a stationary wave is
determined by the system
\begin{equation}
%\aligned
\begin{array}{@{\extracolsep{-.6pc}}rl}
&S = S_- = S_+,\\[5pt]
&p = p(\rho,S_-),\\[5pt]
&[a\rho u] =  0,\\[5pt]
& \Big[\dfrac{u^2}{2} + h(\rho,S_-)\Big] = 0.
%\endaligned
\end{array}
\label{4.5}
\end{equation}
Here, $h  =  \eps + pv$ is the specific enthalpy which satisfies
\begin{equation}
\frac{\pt}{\pt\rho}h(\rho,S_-) \, = \, v
\frac{\pt}{\pt\rho}p(\rho,S_-). \label{4.6}
\end{equation}
To solve the system of equations \rf{4.5} for $\rho_+ = \rho$, we find the roots of the
equation
\begin{equation}
\Phi(\rho) := (u_-^2 + 2 h(\rho_-,S_-)) \rho^2 - 2 \rho^2 h(\rho,S_-) \, = \, \Big(\frac{a_-u_-\rho_-}{a_+}\Big)^2.
\label{4.7}
\end{equation}
A basic calculation implies that the equation \rf{4.7} has a root if and only if
\be
a_+ \, \ge \, a_{\min}(U_-) := \frac{a_-u_-\rho_-}{\sqrt{\Phi(\rho_{\max})}},
\label{4.8}
\ee
where $\rho_{\text{max}}$ is the (unique) value such that ${d\Phi(\rho_{\text{max}})\over d\rho}=0$.
 Moreover, in this case,  \rf{4.7} has two roots $\varphi_1(U_-,a_+) \le
\varphi_2(U_-,a_+)$, which coincide if and only if $a_+  =
a_{\min}(U_-)$.

Next, we address the question of how to select the right states between $\varphi_1(U_-,a_+) $ and $\varphi_2(U_-,a_+)$.  
  First, given a state $U_-$ we observe that the
last equation of \rf{4.5} also determines a {\it stationary curve} $u = u(\rho)$ in
the plane $(\rho, u)$. Hence, the third equation
of \rf{4.5} implies that the component $a$ can be expressed as a
function $a = a(\rho )$ of the variable
 $\rho$ along this curve. We then postulate the following admissibility criterion  (see
\cite{IsaacsonTemple92, LeFlochThanh}).

\bigskip

\textsc{Monotonicity Criterion.}  {\it
  Along the stationary curve in the $(\rho, u)$-plane and between
the left- and right-hand states of any stationary wave, the component
$a$ determined by {\rm \rf{4.5}} and expressed  as a function of $\rho$ must be monotone in $\rho$. }

\bigskip

It is checked in \cite{LeFlochThanh} that: 

\begin{lemma}
\label{Lem4.1}  
The Monotonicity Criterion is equivalent to the condition
that any stationary wave remains in (the closure of) a single phase.
\end{lemma}

\subsection{Definition of the well-balanced scheme}

The numerical scheme considered in the present paper is defined as follows.
The mesh-size is chosen to be uniform, i.e., $x_{j+1} - x_j = \Delta x = l$, and we introduce
standard notation
\be
\aligned
& \frac{1}{\lambda} \geq \max_{j,n} \{|u_j^n| + \sqrt{2 p_\rho(\rho_j^n,S_j^n)}\},
\\
& \Delta t = \lambda\Delta x,\\
&U := (\rho, \rho u, \rho e),\quad f(U) := (\rho u, (\rho u^2 + p), u(\rho e+p)),\\
 U_j^{n+1} & =  \, U_j^n - \lambda (g^\text{LF}(U_{j}^n, U_{j+1,-}^n) - g^\text{LF}(U_{j-1, +}^n, U_{j}^n)),\\
&  = \dfrac{U_{j-1, +}^n + U_{j+1,-}^n}{2} + \frac{\lambda}{2}(f(U_{j-1,+}^n) - f(U_{j+1, -}^n)),\\
\endaligned
 \label{4.9}
\ee
where $g^\text{LF}(U,V)$ is the Lax-Friedrichs numerical flux:
\be
 g^\text{LF}(U,V) =  \frac{1}{2}(f(U) + f(V)) -  \frac{1}{2\lambda}(V - U).
 \label{4.10}
\ee
The description of the states
\be
U_{j+1,-}^n = (\rho, \rho u, \rho e)_{j+1,-}^n,\quad
U_{j-1,+}^n = (\rho, \rho u, \rho e)_{j-1,+}^n
 \label{4.11}
\ee will be given shortly.

\begin{remark}
 Although the mesh-size is chosen to be uniform, the arguments below still hold for non-uniform meshes. 
 For example, one may take $\Delta t = \lambda \inf_{i\in \ZZ}\{|x_{i+1}-x_i|\}$, 
 provided the infimum is not zero.
\end{remark}

In the scheme \rf{4.9}, the states
$$
U_{j+1,-}^n = (\rho, \rho u, \rho e)_{j+1,-}^n,\quad
U_{j-1,+}^n = (\rho, \rho u, \rho e)_{j-1,+}^n
$$
 are defined as follows.

First, we observe that the entropy is constant
across each stationary jump, and we set $\eps_{j+1,-}^n = \eps(\rho_{j+1,-}^n,S_{j+1}^n),
e_{j+1,-}^n = \eps_{j+1,-}^n + (u_{j+1,-}^n)^2/2$, and so on. Then, we
determine
$\rho_{j+1,-}^n, u_{j+1,-}^n$ from the equations
\be
\aligned
& a_{j+1}^n \rho_{j+1}^n u_{j+1}^n = a_j^n \rho_{j+1,-}^n u_{j+1,-}^n,\\
&\dfrac{({u_{j+1}^n})^2}{2} + h(\rho_{j+1}^n, S_{j+1}^n) = \dfrac{({u_{j+1,-}^n})^2}{2} + h(\rho_{j+1,-}^n, S_{j+1}^n),\\
\endaligned
\label{4.12}
\ee
where $h$ denotes the specific enthalpy defined by $dh = TdS + vdp$,
and is expressed as a function $h=h(\rho, S)$ of the density and the specific entropy. In fact, we have
$$
h_\rho(\rho,S) \, = \,  \frac{p_\rho(\rho,S)}{\rho},
$$
for any fixed $S$, and the stationary jump always remains in a given region $G_i$.

 Similarly, we determine $\rho_{j-1,+}^n, u_{j-1,+}^n$ from the equations
\be
\aligned
& a_{j-1}^n \rho_{j-1}^n u_{j-1}^n = a_j^n \rho_{j-1,+}^n u_{j-1,+}^n,\\
&\frac{({u_{j-1}^n})^2}{2} + h(\rho_{j-1}^n, S_{j-1}^n) = \frac{({u_{j-1,+}^n})^2}{2} + h(\rho_{j-1,+}^n, S_{j-1}^n),\\
\endaligned
\label{4.13}
\ee
by requiring that the stationary jump always remains in a given region $G_i$.
This definition selects the state in the same phase. Observe from \eqref{4.13} that
we have two states, one belongs to the same phase (or region) and one belongs to another phase.

\subsection{Preservation of steady states}

Our well-balanced scheme is stable and maintains the equilibrium (steady) states \cite{KroenerThanh2}. 
Considering a stationary wave, then one has  
$$
\begin{array}{@{\extracolsep{-.6pc}}rl}
& a_{j+1}^n \rho_{j+1}^n u_{j+1}^n = a_j^n \rho_j^n u_j^n,\\[6pt]
& \dfrac{({u_{j+1}^n})^2}{2} + h(\rho_{j+1}^n) = \frac{({u_j^n})^2}{2} + h(\rho_j^n).
\end{array}
$$
It follows from \rf{4.12} and \rf{4.14} that
$$
\begin{array}{@{\extracolsep{-.6pc}}rl}
&\rho_{j+1,-}^n = \rho_j^n,\quad u_{j+1,-}^n = u_j^n,\\[5pt]
&\rho_{j-1,+}^n = \rho_j^n,\quad u_{j-1,+}^n = u_j^n,
\end{array}
$$
i.e.
$$
U_{j+1,-}^n = U_j^n, \quad U_{j-1,+}^n = U_j^n,
$$
and, thus
$$
U_j^{n+1} = U_j^n.
$$
This means that the scheme generates the true steady states.

On the other hand, in the special case that $a$ is a constant, then
$$
U_{j+1,-}^n = U_{j+1}^n, \quad U_{j-1,+}^n = U_{j-1}^n,
$$
and our scheme reduces to a rather standard scheme.

Let us next illustrate via some numerical experiments 
the property that the scheme \rf{4.9} preserves steady states.
Consider the Riemann problem for \rf{1.1} where we take $\gamma=1.4$ and the Riemann data are such that 
the left-hand side state is
$$
U_L=(\rho_L, u_L, p_L, a_L) = (2, 0.5, 2^\gamma, 1),
$$
and the right-hand side state
$$
U_R=(\rho_R,u_R,p_R,a_R)=(2.080717229626240,  0.320402338758170, \rho_R^\gamma, 1.5),
$$
which is chosen so that these two states are connected by a stationary wave: 
\be
 U(x,t)=\left\{\begin{array}{ll} U_L,\quad&\text{if } x<0\\
 U_R,\quad&\text{if } x>0.
 \end{array}\right.
\label{4.14}
\ee

We compare here the approximations of the stationary wave \rf{4.14} given by 
our scheme \rf{4.9} and a standard scheme based on a direct discretization of the right-hand side 
(using central difference, forward difference, or backward difference). With the latter, 
stationary waves are not preserved. For definiteness, we choose the Lax-Friedrichs (LF) scheme.

\begin{figure} 
  \includegraphics[width=10truecm]{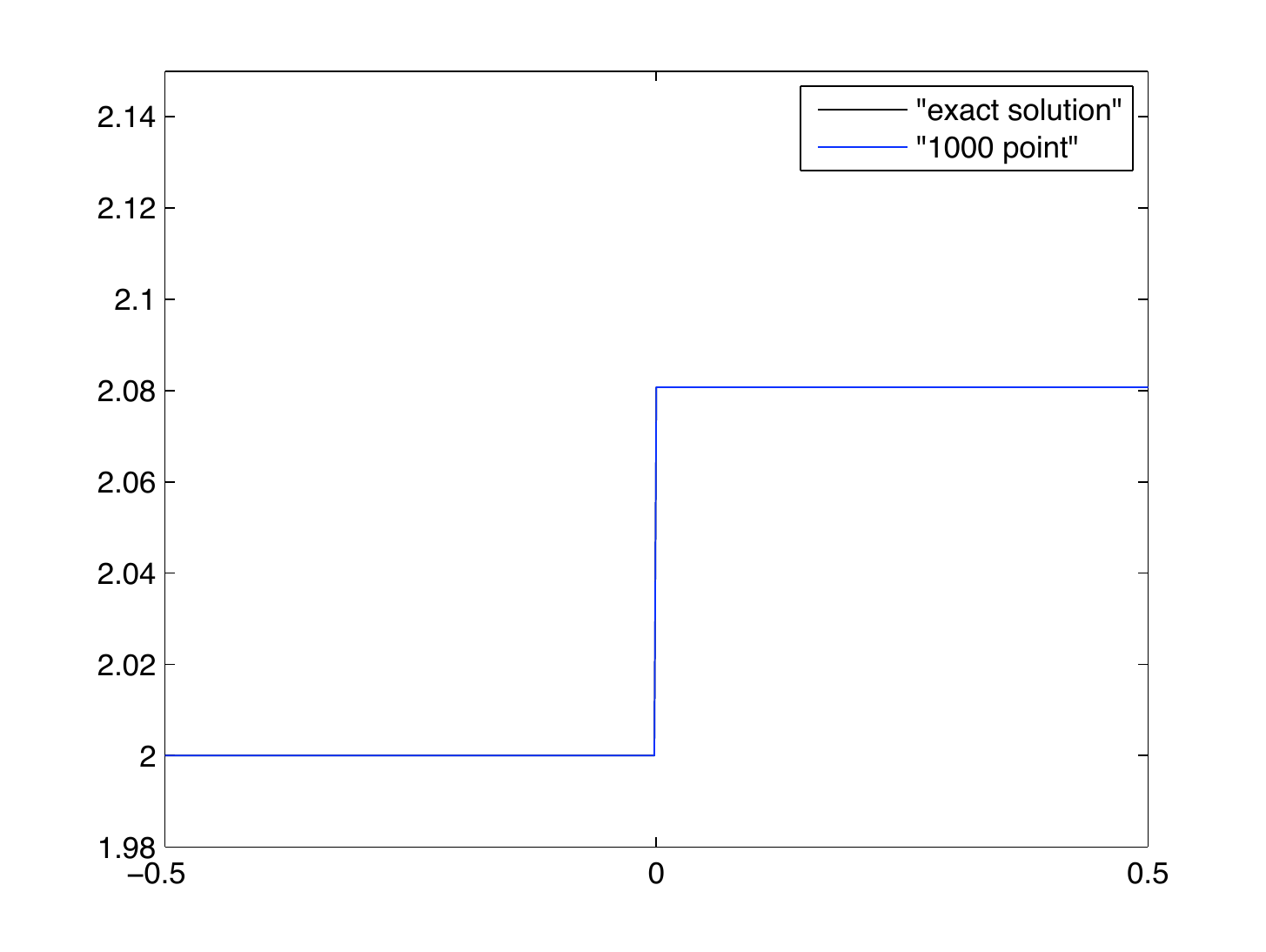}\\
  \caption{Density of the stationary wave \rf{4.14} -- proposed scheme with 1000 mesh points}\label{Fig41}
\end{figure}

\begin{figure} 
  \includegraphics[width=10truecm]{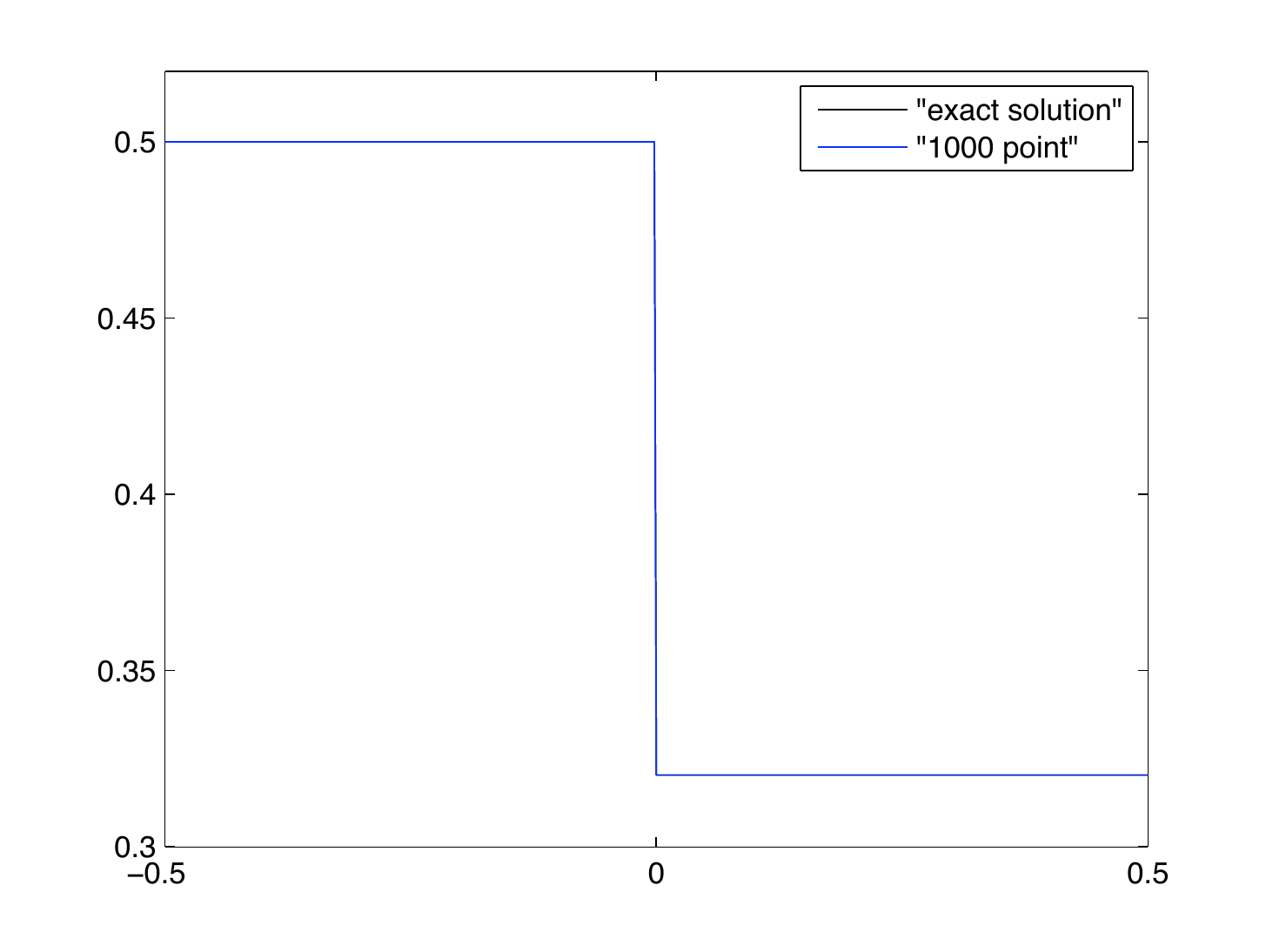}\\
  \caption{Velocity of the stationary wave \rf{4.14} -- proposed scheme with 1000 mesh points}\label{Fig42}
\end{figure}

\begin{figure} 
  \includegraphics[width=10truecm]{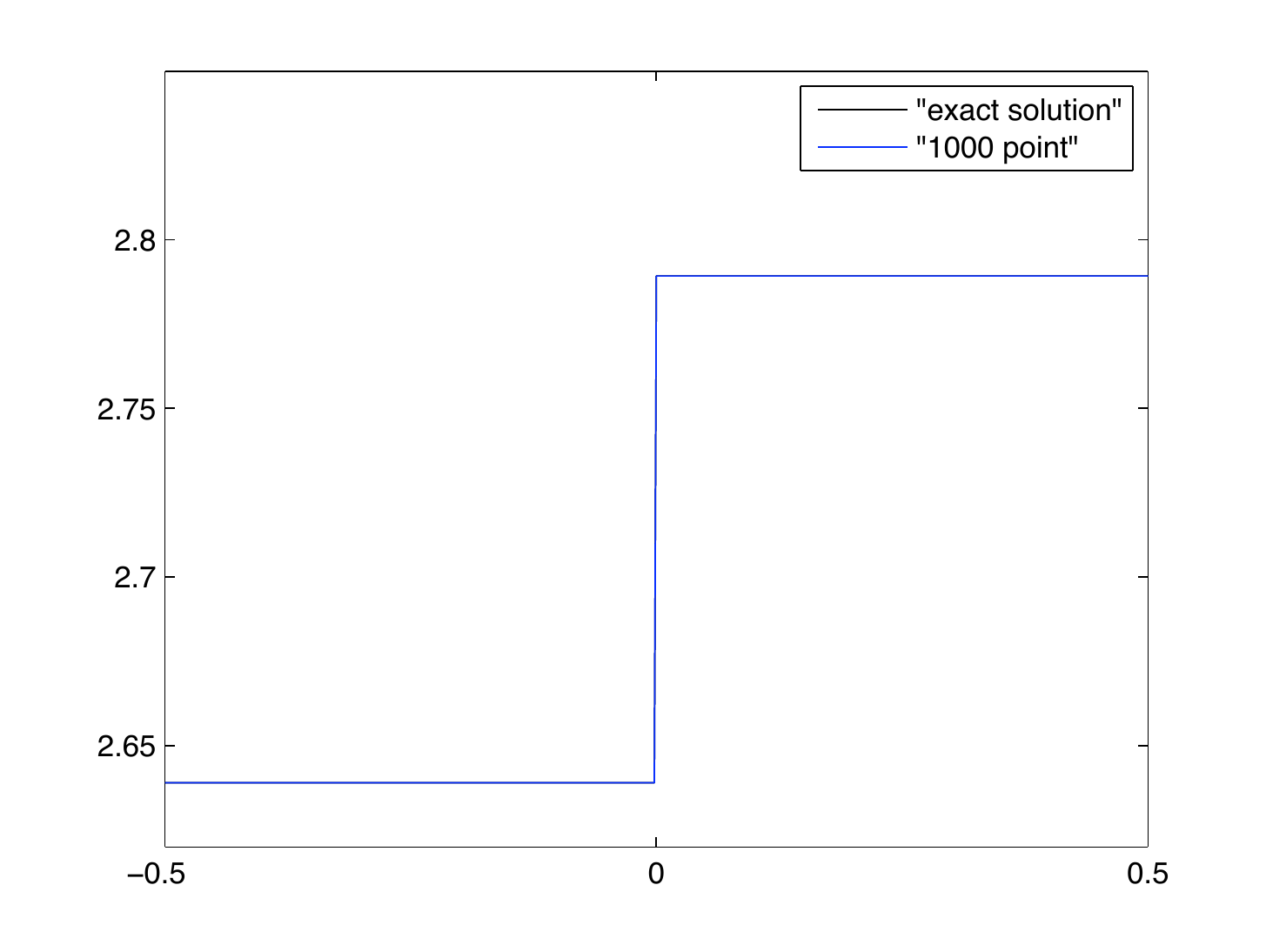}\\
  \caption{Pressure of the stationary wave \rf{4.14} -- proposed scheme with 1000 mesh points}\label{Fig43}
\end{figure}

First, we note that the scheme \rf{4.9} with $1000$ mesh points gives the approximate density, velocity, 
and pressure which are virtually the ones of the exact stationary wave, 
respectively; see Figure \ref{Fig41},  Figure \ref{Fig42}, and Figure \ref{Fig43}.

\begin{figure} 
  \includegraphics[width=10truecm]{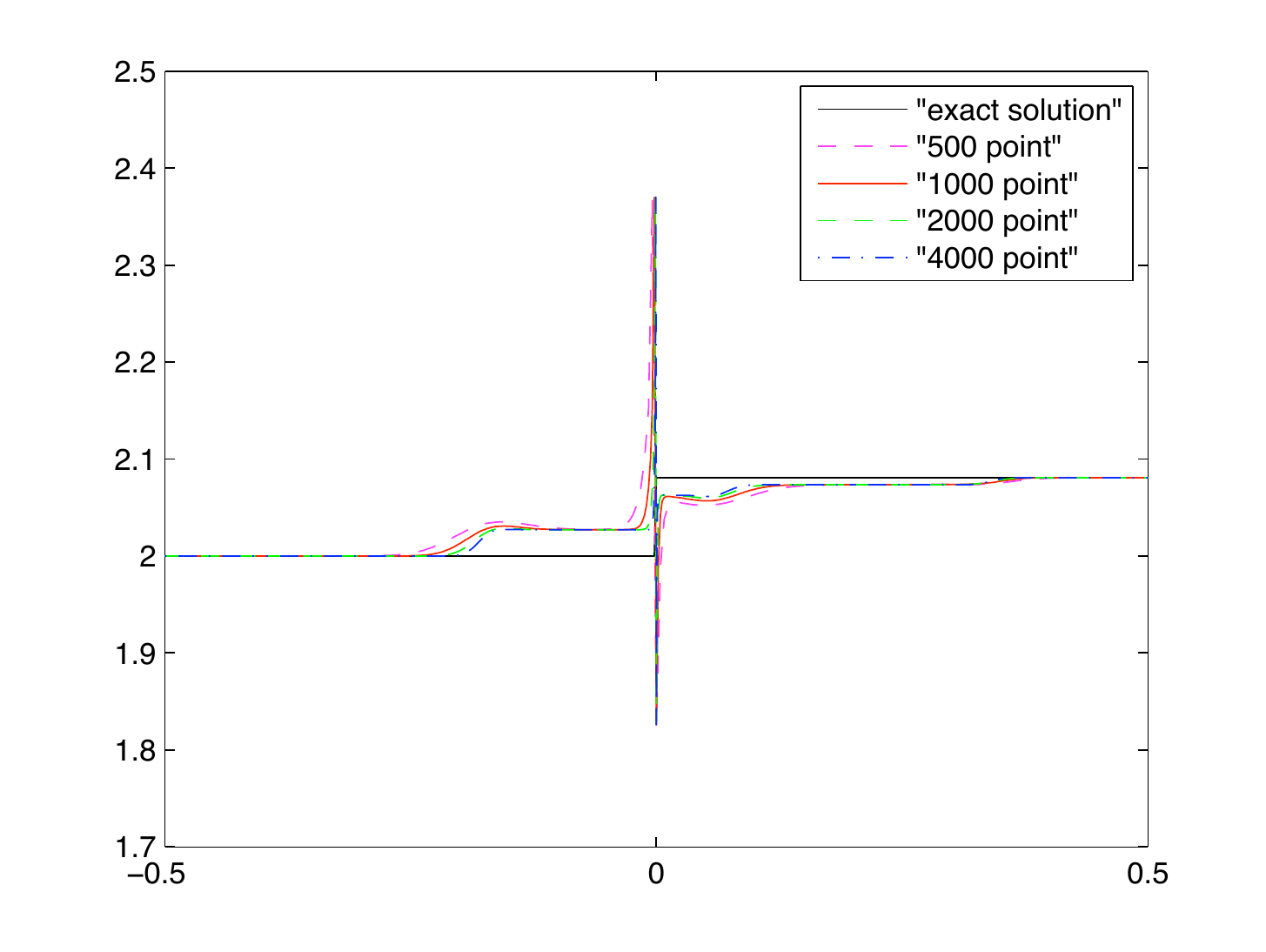}\\
  \caption{Density of the stationary wave \rf{4.14} -- LF 
  scheme with 500, 1000, 2000, and 4000 mesh points}\label{Fig44}
\end{figure}

\begin{figure} 
  \includegraphics[width=10truecm]{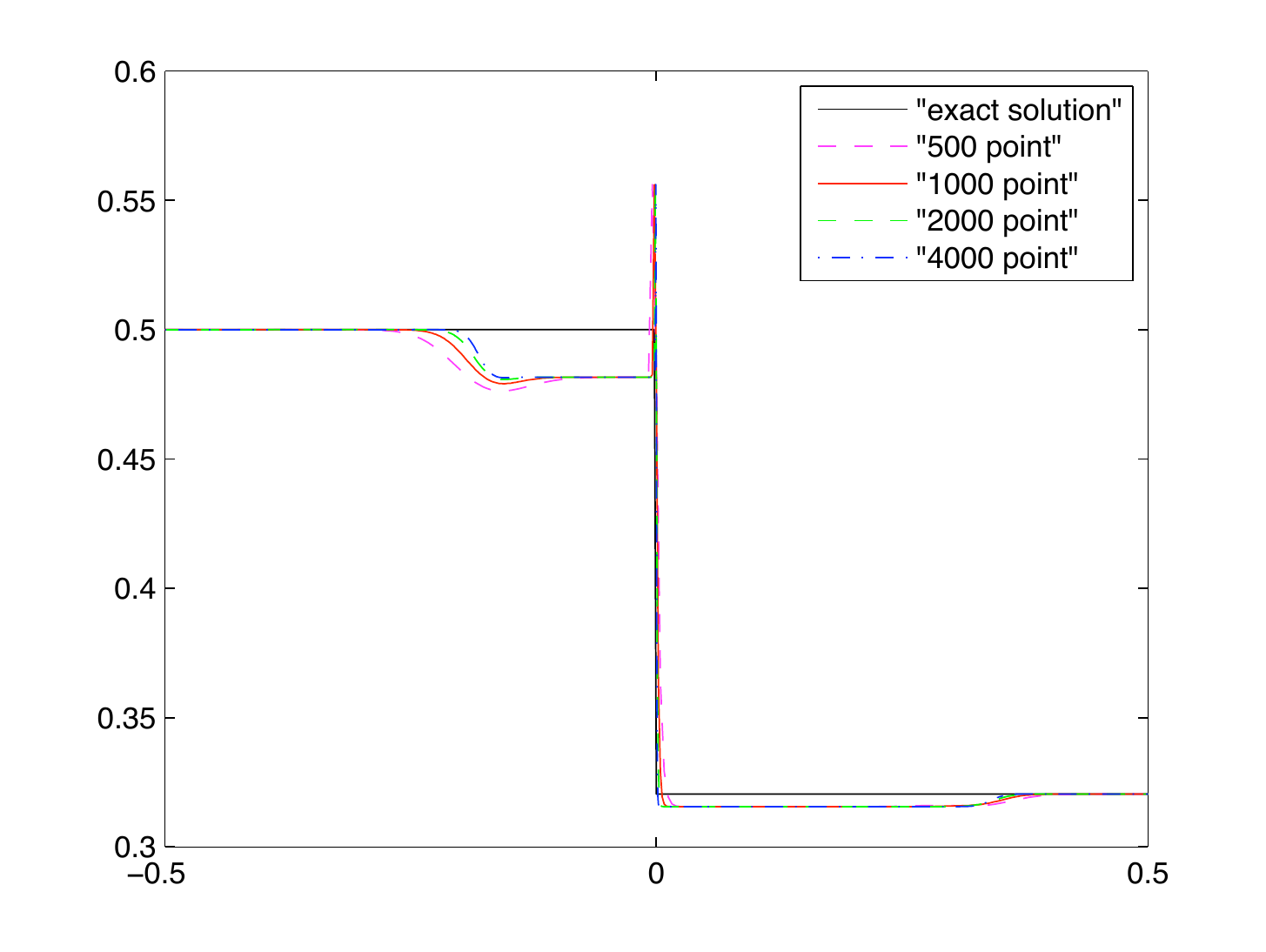}\\
  \caption{Velocity of the stationary wave \rf{4.14} -- LF 
  scheme with 500, 1000, 2000, and 4000 mesh points}\label{Fig45}
\end{figure}

\begin{figure} 
  \includegraphics[width=10truecm]{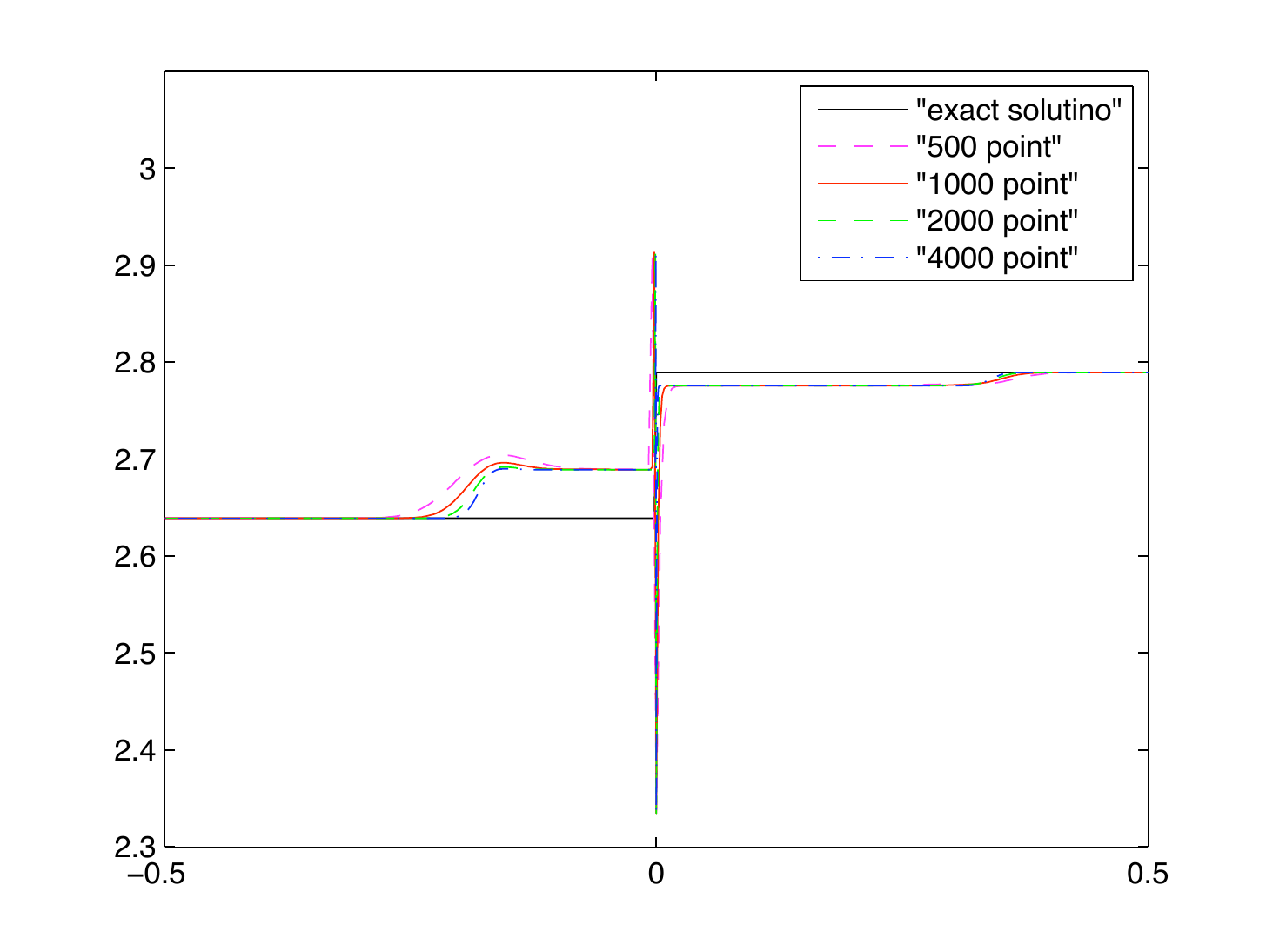}\\
  \caption{Pressure of the stationary wave \rf{4.14} -- LF scheme with 500, 1000, 2000, and 4000 mesh points}\label{Fig46}
\end{figure}

 Second, with the Lax-Friedrichs scheme, the density, velocity, and pressure corresponding to the stationary wave 
 under consideration are not well approximated; see Figure \ref{Fig44},  Figure \ref{Fig45}, Figure \ref{Fig46}. 
 The solution exhibits large spikes near the discontinuity since 
 the equilibrium states are forced out of their equilibrium positions and generate new waves.

%-------------------------------------------------------------------------------------------------------------------

\subsection{Positivity and minimum entropy principle}

Consider stiffened gases, where the local sound speed
$c:=\sqrt{p_\rho(\rho,S)}$ is real. We show that if the initial density if non-negative, then the density of the approximate solution is non-negative. We also 
establish the minimum entropy principle for approximate solutions.

\begin{theorem}\label{theo51}
The scheme  \rf{4.9} satisfies the following : 

1. Positivity property:
 If $\rho_j^0\ge 0$ for all $j$, then $\rho_j^n\ge 0$ for all $j$ and $n$.

2. Minimum entropy principle: for all $j$ and $n$ 
$$
S_j^{n+1} \, \ge \, \min \{S_{j-1}^{n},  S_{j+1}^{n}\}.
$$

\end{theorem}

\begin{proof} 1. We need only show that for any given integer $n$, if $\rho_j^n\ge 0$ for $j$,
 then $\rho_j^{n+1}\ge 0$ for all $j$. Indeed, the scheme \rf{4.9} provides the density
\be
\aligned
 \rho_j^{n+1} & = \frac{\rho_{j-1, +}^n + \rho_{j+1,-}^n}{2}
  + \frac{\lambda}{2}(\rho_{j-1,+}^n u_{j-1,+}^n - \rho_{j+1,-}^n u_{j+1,-}^n)\\
&\ge \frac{\rho_{j-1, +}^n + \rho_{j+1,-}^n}{2}
  - \frac{\lambda}{2} \max\{ |u_{j-1,+}^n|, |u_{j+1,-}^n| \}(\rho_{j-1,+}^n + \rho_{j+1,-}^n)\\
&\ge \frac{\rho_{j-1, +}^n + \rho_{j+1,-}^n}{2} \, \big(1 - \lambda \max\big( |u_{j-1,+}^n|, |u_{j+1,-}^n| \big)\big).
\\
\endaligned
\label{4.15}
\ee
Since stationary waves provided by \rf{4.12} and
\rf{4.14} always connect states with non-negative densities, it
follows from the inequality \rf{4.15} that
\be
 \rho_j^{n+1} > 0 \quad\text{whenever}\quad 1 -  \lambda \max\{ |u_{j-1,+}^n|, |u_{j+1,-}^n| \} > 0.
\label{4.16}
\ee
Besides, the hypothesis implies that
the function $h$ in \rf{4.13} is concave in $\rho$:
$$
 h_{\rho\rho}(\rho,S) =  (\gamma-2)\frac{p_\rho(\rho,S)}{\rho^2} \, \le \, 0.
$$
It follows  from \rf{4.12} and the concavity of the function $h$ that
$$
\aligned
 {1 \over 2} ({u_{j+1,-}^n})^2 - {1 \over 2}({u_{j+1}^n})^2
&= -( h(\rho_{j+1,-}^n, S_{j+1}^n) - h(\rho_{j+1}^n, S_{j+1}^n)) \\
& \le - h_\rho(\rho_{j+1}^n, S_{j+1}^n)(\rho_{j+1,-}^n - \rho_{j+1}^n)\\
& \le - \frac{p_\rho(\rho_{j+1}^n, S_{j+1}^n)}{\rho_{j+1}^n}(\rho_{j+1,-}^n - \rho_{j+1}^n)\\
& \le p_\rho(\rho_{j+1}^n, S_{j+1}^n).
\endaligned
$$
Therefore, we get 
\be 
|u_{j+1,-}^n| \le
\sqrt{(u_{j+1}^n)^2 + 2p_\rho(\rho_{j+1}^n, S_{j+1}^n)} \le
|u_{j+1}^n| + \sqrt{2p_\rho(\rho_{j+1}^n, S_{j+1}^n)}.
\label{5.14} \ee Similarly, \be |u_{j-1,+}^n| \le |u_{j-1}^n| +
\sqrt{2p_\rho(\rho_{j-1}^n, S_{j-1}^n)}. 
\label{5.15} 
\ee
From
\rf{4.9}, \rf{4.16}, \rf{5.14}, and \rf{5.15}, we conclude that
$$
 \rho_j^{n+1} > 0 \quad \text{for all} j,n,
$$
which establishes the first statement of the theorem.

\

2. 
The proof of the second statement is based on the following classical result: 
assume that $\UU$ is a strictly convex function in $\RR^N$, and that there exists a function $\FF$
and a vector-valued map $f$ such that
$D\FF = D\UU\,  Df$.
If $U$ is a vector defined by
\be
U = \frac{V+W}{2} + \frac{\lambda}{2}(f(V)-f(W)),
\label{5.19}
\ee
then
\be
\UU(U) \le  \frac{\UU(V)+\UU(W)}{2} + \frac{\lambda}{2}(\FF(V)-\FF(W)).
\label{5.20}
\ee

Now, comparing \rf{4.9} and \rf{5.19}, we deduce from \rf{5.20} that the scheme \rf{4.9}
satisfies the numerical entropy inequality
$$
 \UU(U_j^{n+1})  \le \frac{\UU(U_{j-1, +}^n) + \UU(U_{j+1,-}^n)}{2} + \frac{\lambda}{2}(\FF(U_{j-1,+}^n) - \FF(U_{j+1, -}^n)),
$$
for any entropy pair of the form \rf{4.1}.
Thus, we have
$$
\aligned
  a_{j}^{n+1} \rho_{j}^{n+1}g(S_{j}^{n+1})  \le  & \frac{1}{2} \Big( 
  a_j^n \rho_{j-1,+}^ng(S_{j-1}^n)
  + a_j^n \rho_{j+1,-}^ng(S_{j+1}^n)\Big)\\
  &+ \frac{\lambda}{2}(a_j^n \rho_{j-1,+}^n u_{j-1,+}^n g(S_{j-1}^n)
   -a_j^n \rho_{j+1,-}^n u_{j+1,-}^n g(S_{j+1}^n)),
   \endaligned
$$
or
$$
  \frac{a_{j}^{n+1}}{a_j^n}  \rho_{j}^{n+1}g(S_{j}^{n+1})
 \le
 \frac{1}{2} \rho_{j-1,+}^n(1   + \lambda u_{j-1,+}^n)g(S_{j-1}^n)
 +
 \frac{1}{2} \rho_{j+1,-}^n(1 - \lambda u_{j+1,-}^n)g(S_{j+1}^n).
$$

In view of Theorem \ref{theo51}, $\rho_{j}^{n+1}$ is
non-negative. Taking $g(S) =  (S_0-S)^p, p>1$, where $S_0$ is some
constant such that $S_0-S>0$, and recalling Lemma~\ref{lem34}
we obtain
$$
\aligned
  \frac{a_{j}^{n+1}}{a_j^n}  \rho_{j}^{n+1}(S_0-S_{j}^{n+1})^p
 & \le
 \frac{1}{2} \rho_{j-1,+}^n(1   + \lambda
 u_{j-1,+}^n)(S_0-S_{j-1}^n)^p
  \\
  &+
 \frac{1}{2} \rho_{j+1,-}^n(1 - \lambda u_{j+1,-}^n)(S_0-S_{j+1}^n)^p.
 \endaligned
 $$
Thus, we get
$$
\aligned
 \Big( \frac{a_{j}^{n+1}}{a_j^n}  \rho_{j}^{n+1}\Big)^{1/p} (S_0-S_{j}^{n+1})
 & \le
\Big(\frac{1}{2} \rho_{j-1,+}^n(1   + \lambda u_{j-1,+}^n)
(S_0-S_{j-1}^n)^p\\
 &\qquad +
 \frac{1}{2} \rho_{j+1,-}^n(1 - \lambda u_{j+1,-}^n)(S_0-S_{j+1}^n)^p\Big)^{1/p}\\
 & \le
 \Big(\frac{1}{2} \rho_{j-1,+}^n (1   + \lambda u_{j-1,+}^n)
  +
\frac{1}{2} \rho_{j+1,-}^n(1 - \lambda u_{j+1,-}^n)\Big)^{1/p}\\
&\qquad\times
\max\{S_0-S_{j-1}^n, S_0-S_{j+1}^n\}.\\
\endaligned
 $$
Letting $p\to +\infty$ in this inequality we obtain
$$
S_0-S_{j}^{n+1}\le \max\{S_0-S_{j-1}^n,
S_0-S_{j+1}^n\}=S_0-\min\{S_{j-1}^n, S_{j+1}^n\},
$$
or
$$
S_{j}^{n+1}\ge \min\{S_{j-1}^n, S_{j+1}^n\},
$$
which completes the proof.
\end{proof}

%================================================================================================================


\begin{thebibliography}{10}

\bibitem{AndrianovWarnecke}
N.~Andrianov and G.~Warnecke,
\newblock {On the solution to the Riemann problem for the compressible duct
  flow},
\newblock {\em SIAM J. Appl. Math.}, 64(3):878--901, 2004. 

\bibitem{AudusseBouchutBristeauKleinPerthame}
E.~Audusse, F.~Bouchut, M-O. Bristeau, R.~Klein, and B.~Perthame,
\newblock A fast and stable well-balanced scheme with hydrostatic
  reconstruction for shallow water flows,
\newblock {\em SIAM J. Sci. Comp.}, 25(6):2050--2065, 2004. 

\bibitem{BotchorishviliPerthameVasseur}
R.~Botchorishvili, B.~Perthame, and A.~Vasseur,
\newblock {Equilibrium schemes for scalar conservation laws with stiff
  sources},
\newblock {\em Math. Comput.}, 72:131--157, 2003.

\bibitem{BotchorishviliPironneau}
R.~Botchorishvili and O.~Pironneau,
\newblock Finite volume schemes with equilibrium type discretization of source
  terms for scalar conservation laws,
\newblock {\em J. Comput. Phys.}, 187:391--427, 2003.

\bibitem{Bouchut}
F. Bouchut,
\newblock Nonlinear stability of finite volume methods for hyperbolic conservation laws,
and well-balanced schemes for sources,
Frontiers in Mathematics series, Birkh\"auser, 2004.

\bibitem{CourantFriedrichsbook}
R.~Courant and K.O. Friedrichs,
\newblock {\em Supersonic flow and shock waves},
\newblock John Wiley, New York, 1948.

\bibitem{DalMasoLeFlochMurat}
G.~Dal~Maso, P.G. LeFloch, and F.~Murat,
\newblock Definition and weak stability of nonconservative
products,
\newblock {\em J. Math. Pures Appl.}, 74:483--548, 1995.

\bibitem{GoatinLeFloch}
P.~Goatin and P.G. LeFloch,
\newblock {The Riemann problem for a class of resonant nonlinear systems of
  balance laws},
\newblock {\em Ann. Inst. H. Poincar\'e Anal. NonLin\'eaire}, 21:881--902, 2004.

\bibitem{Gosse00}
L.~Gosse,
\newblock A well-balanced flux-vector splitting scheme designed for hyperbolic
  systems of conservation laws with source terms,
\newblock {\em Comp. Math. Appl.}, 39:135--159, 2000.

\bibitem{GreenbergLeroux}
J.M. Greenberg and A.Y. Leroux,
\newblock A well-balanced scheme for the numerical processing of source terms
  in hyperbolic equations,
\newblock {\em SIAM J. Numer. Anal.}, 33:1--16, 1996.

\bibitem{HartenLaxLevermoreMorokoff98}
A.~Harten, P.D. Lax, C.D. Levermore, and W.J. Morokoff,
\newblock {Convex entropies and hyperbolicity for general Euler
equations},
\newblock {\em SIAM J. Numer. Anal.}, 35(6):2117--2127, 1998.

\bibitem{IsaacsonTemple92}
E.~Isaacson and B.~Temple,
\newblock Nonlinear resonance in systems of conservation laws,
\newblock {\em SIAM J. Appl. Math.}, 52:1260--1278, 1992.

\bibitem{IsaacsonTemple95}
E.~Isaacson and B.~Temple,
\newblock Convergence of the $2\times 2$ Godunov method for a general resonant
  nonlinear balance law,
\newblock {\em SIAM J. Appl. Math.}, 55:625--640, 1995.

\bibitem{KroenerThanh1}
D.~Kr\"oner and M.D. Thanh,
\newblock {On the Model of Compressible Flows in a Nozzle: Mathematical
  Analysis and Numerical Methods},
\newblock {\em Proc. Tenth. Intern. Conf. "Hyperbolic Problem: Theory, Numerics,
and Applications", Osaka (2004)},
 Yokohama Publishers, 117-124, 2006.

\bibitem{KroenerThanh2}
D.~Kr\"oner and M.D. Thanh,
\newblock {Numerical solutions to compressible flows in a nozzle with variable
  cross-section},
\newblock {\em SIAM J. Numer. Anal.}, 43(2): 796-824, 2006.

\bibitem{Lax71}
P.D. Lax,
\newblock {Shock waves and entropy, in: E.H. Zarantonello, Ed.},
\newblock {\em Contributions to Nonlinear Functional Analysis},  603--634,
  1971.

\bibitem{LeFloch88}
P.G. LeFloch,
\newblock {Entropy weak solutions to nonlinear hyperbolic systems under
  nonconservative form},
\newblock {\em Com. Partial. Diff. Eqs.}, 13(6):669--727, 1988.

\bibitem{LeFlochbook}
P.G. LeFloch,
\newblock {\em Hyperbolic systems of conservation laws:
The theory of classical and non-classical shock waves}.
\newblock Lectures
in Mathematics, ETH Z\"urich, Birk\"auser, 2002.

\bibitem{LeFloch89}
P.G. LeFloch,
\newblock {Shock waves for nonlinear hyperbolic systems in nonconservative
  form},
\newblock {\em Institute for Math. and its Appl., Minneapolis, Preprint}, 593,
  1989.

\bibitem{LeFloch04}
P.G. LeFloch,
\newblock{Graph solutions of nonlinear hyperbolic systems,}
\newblock{\em J. Hyper. Diff. Equa.}, 1:243--289, 2004.

\bibitem{LeFlochLiu}
 P.G. LeFloch and T.-P. Liu,
 \newblock{Existence theory for nonlinear hyperbolic systems in nonconservative form,}
 \newblock{\em Forum Math.}, 5:261--280, 1993.

\bibitem{LeFlochThanh}
P.G. LeFloch and M.D. Thanh,
\newblock {The Riemann problem for fluid flows in a nozzle with discontinuous
  cross-section},
\newblock {\em Comm. Math. Sci.}, 1(4):763--797, 2003.

\bibitem{LeFlochThanh2}
P.G. LeFloch and M.D. Thanh,
\newblock {The Riemann problem for the shallow water equations with
discontinuous topography}, 
\newblock {\em Comm. Math. Sci.}, (to appear).

\bibitem{MarchesinPaes-Leme}
D.~Marchesin and P.J. Paes-Leme,
\newblock {A Riemann problem in gas dynamics with bifurcation. Hyperbolic
  partial differential equations III},
\newblock {\em Comput. Math. Appl. (Part A)}, 12:433--455, 1986.

\bibitem{Tadmor84}
E.~Tadmor,
\newblock {Skew selfadjoint form for systems of conservation
laws},
\newblock {\em J. Math. Anal. Appl.}, 103:428--442, 1984.

\bibitem{Tadmor86}
E.~Tadmor,
\newblock {A minimum entropy principle in the gas dynamics
equations},
\newblock {\em Appl. Numer. Math.}, 2:211--219, 1986.

\end{thebibliography}
\end{document}